\documentclass{article}
\font\tenfrakt=eufm10   \font\tenscript=eusm10
\font\sevenfrakt=eufm7  \font\sevenscript=eusm7
\font\fivefrakt=eufm5   \font\fivescript=eusm5

\newfam\fraktfam
\textfont\fraktfam=\tenfrakt
\scriptfont\fraktfam=\sevenfrakt
\scriptscriptfont\fraktfam=\fivefrakt

\newfam\scriptfam
\textfont\scriptfam=\tenscript
\scriptfont\scriptfam=\sevenscript
\scriptscriptfont\scriptfam=\fivescript

\font\tenbbb=msbm10   \font\tenscript=msbm10
\font\sevenbbb=msbm7  \font\sevenscript=msbm7
\font\fivebbb=msbm5   \font\fivescript=msbm5

\newfam\bbbfam
\textfont\bbbfam=\tenbbb
\scriptfont\bbbfam=\sevenbbb
\scriptscriptfont\bbbfam=\fivebbb
\def\bbb{\fam\bbbfam\tenbbb}

\def\R{{\bbb R}}
\def\C{{\bbb C}}
\def\Z{{\bbb Z}}
\def\N{{\bbb N}}

\newtheorem{theorem}{Theorem}

\newtheorem{lemma}[theorem]{Lemma}
\newtheorem{proposition}[theorem]{Proposition}
\newtheorem{corollary}[theorem]{Corollary}

\title{Cohomological Restrictions on K\protect \"ahler groups}
\author{\sc Azniv  Kasparian   }
\date{      }

\begin{document}
\maketitle

\footnotetext{


   Partially supported by Bulgarian Ministry of Education,
Grant MM 1003 / 2000 .}

\begin{abstract}
 Let $ X $ be a compact K\"ahler manifold with fundamental group $ \pi_{1}(X). $
After introducing the notion of higher Albanese genera $ g_{k}, $ the work
 establishes lower bounds on the number of the relations of $ \pi_{1}(X) $
in terms of the number of the generators, the irregularity, the Albanese dimension,
$ g_{k} $  and etc. The argument relates the cup product maps in the cohomologies of $ X $
and $ \pi _{1}(X). $ It derives some lower bounds on the ranks of these cup products and applies
 Hopf's Theorem, describing $ H_{2}( \pi_{1}(X), {\Z}). $ The same techniques provide lower
  bounds on the Betti numbers of $ X $ and $ \pi _{1}(X) $ within the range of the
  Albanese dimension.
 \end{abstract}

 \section{Statement of the results}

 \hspace{.8mm}

          The abstract groups $ G $ which are isomorphic to the fundamental group
$ \pi_{1}(X) $ of a compact K\"ahler manifold $ X $ are briefly referred to as
 K\"ahler groups. These are always finitely presented.

      The compact complex torus $ Alb(X) = H^{1,0}(X)^{*} /
H_{1}(X,{\Z})_{free} $ is called an Albanese variety of the compact K\"ahler
manifold $ X. $ The  Albanese map $ alb _{X} : X \rightarrow Alb(X), $
$ alb_{X} (x) ( \omega ):= \int _{x_{0}}^{x} \omega $ for $ \omega \in H^{1,0}(X) $
is defined up to a translation, depending on the choice of a base point
$ x_{0} \in X. $ The Albanese dimension of $ X $ is $ a = a(X) :=
 \dim _{\C} alb_{X}(X). $

       The compact K\"ahler manifold $ Y $ is said to be Albanese general if
$ h^{1,0}(Y) >  \dim _{\C} Y = a(Y). $ A surjective holomorphic map
$ f_{k} : X \rightarrow Y_{k} $ of a compact K\"ahler manifold $ X $ onto an Albanese
general manifold $ Y_{k} $ of $ \dim _{\C} Y_{k} = k $  is called an Albanese general
$k$-fibration. It  induces a complex linear embedding
 $ f_{k}^{*} : H^{1,0}(Y_{k}) \rightarrow H^{1,0}(X) $ of the holomorphic $(1,0)$-forms,
  so that $ h^{1,0}(Y_{k}) $ is bounded above by  $  h^{1,0}(X). $ The maximal
 $ h^{1,0}(Y_{k}) $ for Albanese general $k$-fibrations $ f_{k}: X \rightarrow Y_{k} $
  is called $k$-th Albanese genus of $ X $ and denoted by  $ g_{k} = g_{k}(X). $

       The aim of the present note is to establish the following estimates:

\begin{proposition}      \label{GenRel}
Let $ X $ be a compact K\"ahler manifold whose fundamental group
 admits a finite presentation $ \pi_{1}(X) = F / R $ where
 $ F := \langle x_{1}, \ldots , x_{s} \rangle $ is a free group on $ s $ generators,
  $ R_{o} := \langle y_{1}, \ldots , y_{r}  \rangle $ is the subgroup of $ F, $
  generated by the relations $ y_{1}, \ldots , y_{r} \in F $ and $ R $ is the
  normal subgroup of $ F, $ generated by $ R_{o}. $ Suppose that the subgroup
 $ K := ( R_{o} \cap [F,R]) / [ R_{o}, R_{o}] $ of the abelianization
 $ ab R_{o} := R_{o} / [R_{o}, R_{o}] \simeq {\Z}^{r} $ is of $ rk K = k, $
  $ h^{1,0} := \dim _{\C} H^{1,0}(X) $ is the irregularity of $ X, $
  $ a $ is the Albanese dimension and $ g_{k}, $ $ 1 \leq k \leq n $
  are the Albanese genera of $ X. $ Then

  (i) $ r \geq s + k $ for $ h^{1,0} = 0, $ $ a = 0 ; $

  (ii) $ r \geq s + k - 2 h^{1,0} + 1 $ for $ h^{1,0} \geq 1, $ $ a = 1; $

  (iii) $ r \geq s + k - 2 h^{1,0} + \max ( a(a-1), g_{k}(g_{k} -1)  \ \ \vert \ \
   2 \leq k \leq a ) \ \  + \ \  $

    $ \max \left( \frac{a(a-1)}{2}, 2a-1, g_{k} -1 \ \ \vert \ \  2 \leq k \leq a
  \right) $  for $ h^{1,0} \geq g_{1} \geq 2, $ $ a \geq 2; $

  (iv) $ r \geq s + k - 2 h^{1,0} + \max ( 4 h^{1,0} - 6, a(a-1), g_{k}(g_{k}-1)
  \ \ \vert \ \ 2 \leq k \leq a ) \ \ + \ \  $

  $ \max \left( 2 h^{1,0} - 1, \frac{a(a-1)}{2},  g_{k} - 1 \ \ \vert \ \
  2 \leq k \leq a \right) $ for $ h^{1,0} \geq 2, $ $ a \geq 2, $ $ g_{1} = 0. $
  \end{proposition}

            A K\"ahler group $ \pi_{1}(X) $ admits various finite presentations and
 there is no general algorithm for deciding whether two presentations
 determine isomorphic groups. The aforementioned Proposition \ref{GenRel} is not
 expected to perform K\"ahler tests on abstract finitely presented groups. It rather
 studies the influence of some cohomological properties of compact K\"ahler manifolds
 $ X $ on their fundamental groups $ \pi_{1}(X). $ Part of the techniques from the
 proof of Proposition \ref{GenRel} provide also the following

\begin{proposition}   \label{BettiNum}
Let $ X $ be a compact K\"ahler manifold with positive irregularity  $ h^{1,0}, $
 Albanese dimension $ 1 \leq a \leq n = \dim _{\C} X $ and Albanese genera
 $ g_{k}, $ $ 1 \leq k \leq n. $ Then the Betti numbers
 $ b_{m} ( \pi_{1}(X)) := \dim _{\C} H^{m} ( \pi_{1}(X), {\C}) $ and
 $ b_{m}(X) := \dim _{\C} H^{m} ( X, {\C} ) $ are bounded below as follows :

 $$ b_{2i}( \pi_{1}(X)) \geq 2 \sum _{j=0}^{i-1} \mu ^{j,2i-j} + \mu ^{i,i},  \ \
 b_{2i+1} ( \pi_{1}(X)) \geq 2 \sum _{j=0}^{i} \mu ^{j,2i+1-j} \ \ \mbox{  for  }
  3 \leq 2i, 2i+1 \leq a, $$
 $$ b_{2i}(X) \geq 2 \sum _{j=0}^{i-1} \mu ^{j,2i-j} + \mu ^{i,i}, \ \
 b_{2i+1} (X) \geq 2 \sum _{j=0}^{i} \mu ^{j,2i+1-j} \ \ \mbox{  for  }
 3 \leq 2i, 2i+1 \leq a, $$
 $$ b_{2i}(X) \geq 2 \sum _{j=0}^{i-1} \mu ^{n-2i+j,n-j} + \mu ^{n-i,n-i}, \ \
 b_{2i+1}(X) \geq 2 \sum _{j=0}^{i} \mu ^{n-2i-1+j,n-j} \ \ \mbox{  for }
 2n-a \leq 2i, 2i+1 \leq 2n-3, $$
 where
 $$ \mu ^{i,j} := \max \left( \left( \begin{array}{c}
                                                  a \\
                                                  i+j
                                          \end{array} \right),
\left( \begin{array}{c}
          g_{k} - i \\
          j
       \end{array} \right), \delta ^{0}_{g_{1}} \ldots \delta ^{0}_{g_{i+j-1}}
[ (i+j)(h^{1,0} - i - j) + 1] \ \ \vert \ \ g_{k} > 0, i+j \leq k \leq a \right) $$
for $ i \leq j $ and $ \delta ^{0}_{g_{s}} $ standing for Kronecker's delta.
\end{proposition}

       The next section specifies the cases in which the bounds from Proposition \ref{GenRel}
 are stronger than the already known results.  Section 3 collects some properties of
 Albanese dimension and Albanese genera, necessary for deriving lower bounds on the
 ranks of  cup products in $ H^{*} (X,{\C}). $ Section 4 relates cup
 products in  group cohomologies $ H^{*} ( \pi_{1}(X), {\C} ) $ with the
 corresponding cup products in de Rham cohomologies $ H^{*} (X,{\C}). $ Section 5
 justifies that $ \mu ^{i,m-i} $ from Proposition \ref{BettiNum} are lower bounds on
 the ranks of the cup products
 $ \zeta ^{i,m-i}_{X} : \wedge ^{i} H^{1,0}(X) \otimes _{\C} \wedge ^{m-i} H^{0,1}(X)
 \rightarrow H^{m}(X, {\C}). $
 The last section 6 recalls Hopf's Theorem on the second homologies of a group
 and concludes the proofs of Propositions \ref{GenRel} and \ref{BettiNum}.

 {\bf Acknowledgements:} The author is extremely grateful to Tony Pantev for bringing
 to her attention the article \cite{A} of Amor\'os and for the useful advices, comments
 and conversations. She apologizes for declaring in the previously circulated version
 a wrong counterexample to a theorem of Remmert and Van de Ven ,and announcing, in this
 way, a nonexisting error in Amor\'os' work \cite{A}. The author thanks Prof. Amor\'os
 for pointing out the aforementioned mistake and explaining her that both Remmert and
 Van de Ven's Theorem and Amor\'os' results \cite{A} are completely accurate.

\section{ Comparison with previous related works}

\hspace{.8mm}

       Prior to Proposition \ref{GenRel} are known the following estimates
       among the number of the generators and relations of a K\"ahler group.

\begin{theorem}    \label{GreenLazarsfeld}
{\rm (Green and Lazarsfeld \cite{GL})} Let $ X $ be a compact K\"ahler manifold
whose fundamental group $ \pi_{1}(X) $ admits a presentation with  $ s $ generators and
 $ r $ relations.

       (i) If the Albanese genus $ g_{1} = 0 $ then $ r \geq s - 3. $

       (ii) If the Albanese dimension $ a \geq 2 $ then $ r \geq s - 1. $
\end{theorem}

For $ s - r \geq 2, $ Green and Lazarsfeld show that the entire character variety
$ \widehat{\pi_{1}(X)} = Hom ( \pi_{1}(X), {\C}^{*}) \simeq H^{1}(X,{\C})  $ of
 $ \pi_{1}(X) $ is contained in the special locus
$ S^{1}(X):= \{ L \in Pic ^{o}(X) \vert H^{1}(X,L) \neq 0 \} $ of the topologically
 trivial line bundles on $ X, $ parametrized by $ Pic ^{o}(X) \subset
  H^{1}(X,{\cal O}^{*}_{X}). $ Then there is a surjective holomorphic map
  $ X \rightarrow C $ onto a curve $ C $ of genus $ \geq \frac{s-r}{2} $ and the
 Albanese image of $ X $ is a curve. That violates the assumption of the second part.
 The first part is contradicted by $ s - r \geq 4, $ as far as $ g_{1} = 0 $
  signifies the nonexistence of surjective holomorphic maps $ X \rightarrow C $
  onto curves $ C $ of genus $ \geq 2. $

 \begin{theorem}   \label{Amoros}
{\rm (Amor\'os \cite{A})} Let $ X $ be a compact K\"ahler manifold with $ g_{1} = 0, $
whose fundamental group $ \pi_{1}(X) $ admits a presentation with $ s $ generators
and  $ r $ relations. Then

 (i) $ r \geq s $ for $ h^{1,0} = 0; $

 (ii) $ r \geq s - 1 $ for $ h^{1,0} = 2; $

 (iii) $ r \geq s + 4h^{1,0} - 7 $ for $ h^{1,0} \geq 2. $
 \end{theorem}

  Let $ \zeta ^{2}_{X} : \wedge ^{2} H^{1} (X,{\C}) \rightarrow H^{2}(X,{\C}) $ be
 the cup product in de Rham cohomologies. Making use of Sullivan's 1-minimal models,
 Amor\'os identifies $ Ker \zeta ^{2}_{X} $ with $ \left( \pi_{1}(X)_{2} / \pi_{1}(X)_{3}
 \right) \otimes {\R} $ where $ \pi_{1}(X)_{1} := \pi_{1}(X), $ $ \pi_{1}(X)_{i+1}
  := [ \pi_{1}(X)_{i}, \pi_{1}(X)] $ are the components of the lower central series
  of the fundamental group $ \pi_{1}(X). $

         For an arbitrary group $ G, $ let $ J_{G} := \{ \sum _{g} r_{g} g \vert
r_{g} \in {\R}, \sum _{g} r_{g} = 0 \} $ be the augmentation ideal of the group ring
$ {\R} [G]. $ It is well known (cf.\cite{W}) that $ H_{1}(G,{\Z}) \simeq G / [G,G]
\simeq J_{G} / J^{2}_{G}. $ Amor\'os shows in \cite{A} that the ${\R}$-linear map
$ \Delta _{1} : \oplus _{j=1}^{r} {\R} y_{j} \rightarrow J_{F} / J_{F}^{2}, $
$ \Delta _{1} (y_{j}) = y_{j} - 1 + J^{2}_{F} $ has $ Coker \Delta _{1} =
J_{\pi_{1}(X)} / J^{2}_{\pi_{1}(X)}. $ In particular, $ \dim _{\R} Ker \Delta _{1}
 = r - s + 2 h^{1,0}. $ Further, the induced map $ \Delta _{2} : Ker \Delta _{1}
 \rightarrow \wedge ^{2} \left( J_{\pi_{1}(X)} / J^{2}_{\pi_{1}(X)} \right) $ is
  proved to have $ Coker \Delta _{2} \simeq \left( \pi_{1}(X)_{2} / \pi_{1}(X)_{3}
  \right) \otimes {\R}. $ Consequently, $ \dim _{\R} Ker \zeta ^{2}_{X} = \dim _{\R}
  \left( \pi_{1}(X)_{2} / \pi_{1}(X)_{3} \right) \otimes {\R} =
  \left( \begin{array}{c}
            2 h^{1,0}  \\
            2
   \end{array} \right) - \dim _{\R} Ker \Delta _{1} + \dim _{\R} Ker \Delta _{2}
    \geq \left( \begin{array}{c}
                   2 h^{1,0}   \\
                   2
                   \end{array}  \right) - r + s - 2 h^{1,0} $   or
$ rk \zeta ^{2}_{X} \leq r - s + 2 h^{1,0}. $

       On the other hand, Amor\'os  makes use of the following lower bounds on the
 ranks of the cup products:

 \begin{lemma}    \label{Amoros1}
 If $ X $ is a compact K\"ahler manifold with $ g_{1} = 0 $ then

 {\it (i)} $ rk [ \zeta ^{2,0}_{X} : \wedge ^{2} H^{1,0}(X) \rightarrow H^{2}(X,{\C})]
 \geq 2 h^{1,0} - 3; $

 {\it (ii)} $ rk [ \zeta ^{1,1}_{X} : H^{1,0}(X) \otimes _{\C} H^{0,1}(X) \rightarrow
 H^{2}(X,{\C}) ] \geq 2 h^{1,0} - 1. $
 \end{lemma}

        The estimate {\it (i)} is derived from the transversality of the cone
$ {\cal C} ^{2,0} := \{ \omega _{1} \wedge \omega _{2} \vert \omega _{1},
 \omega _{2} \in H^{1,0}(X) \} $  to the kernel of $ \zeta ^{2,0}_{X}, $ i.e.,
  $ {\cal C}^{2,0} \cap Ker \zeta ^{2,0}_{X} = \{ 0 \} $ (cf. also \cite{BPV}.)
   The inequality {\it (ii)}  is a consequence of a theorem of Remmert and Van de Ven
   from \cite{RV}. It  asserts that  a holomorphic map $ \tau : A_{1} \times
   A_{2} \rightarrow B $ of projective algebraic manifolds $ A_{1}, A_{2} $ in a
   complex space $ B $ has $ rk _{\C} \tau \geq \dim _{\C} A_{1} + \dim _{\C} A_{2}, $
   provided $ b_{2}(A_{1}) = b_{2}(A_{2}) = 1, $ $ b_{1} (A_{i}) = 0 $ for some
   $ 1 \leq i \leq 2 $ and $ \tau $ does not factor through a canonical projection
   $ \Pi _{i} : A_{1} \times A_{2} \rightarrow A_{i}. $ This is applied to the
   projectivization $ {\bf P} ( \zeta ^{1,1}_{X}) : {\bf P} ( H^{1,0}(X)) \times
   {\bf P} ( H^{0,1}(X)) \rightarrow {\bf P} ( H^{2} (X, {\C})) $ of the bilinear map
   $ \zeta ^{1,1}_{X} $ with trivial kernel.

  Obviously, Theorem \ref{Amoros} implies  Theorem \ref{GreenLazarsfeld}{\it (i)}.
  One checks straightforward that the inequalities, given by {\it (i)},
   {\it (ii)} with $ h^{1,0} = 1 $ and {\it (iv)} from Proposition \ref{GenRel}
   are stronger than the corresponding estimates from Theorem \ref{Amoros}.
   In the case of $ a \geq 2, $ $ h^{1,0} \geq g_{1} \geq 2 $ with $ h^{1,0} $ comparatively
    large with respect  to $ a $ and $ g_{k}, $ $ 2 \leq k \leq a, $ the bound from
     Proposition   \ref{GenRel} {\it (iii)} may happen to be weaker than the one from
     Theorem   \ref{GreenLazarsfeld} {\it (ii)}.

   \section{Preliminaries on Albanese dimension and Albanese genera}

\hspace{.8mm}

         Some of the bounds on the ranks of cup products, proved in section 5,
 require the characterization the Albanese dimension in terms of cup products.

   \begin{proposition}       \label{AlbDim}
{\rm (Catanese \cite{ABCKT},\cite{C})}  The Albanese dimension $ a:= \dim _{\C}
   alb_{X}(X) $ of a compact K\"ahler manifold $ X $ is the greatest integer with
 $ Im [ \zeta ^{a,a}_{X} : \wedge ^{a} H^{1,0}(X) \otimes _{\C} \wedge ^{a} H^{0,1}(X)
  \rightarrow H^{2a}(X, {\C})] \neq 0 $ or, equivalently, the greatest integer with
  $ Im [ \zeta ^{a,0}_{X} : \wedge ^{a} H^{1,0}(X) \rightarrow H^{a}(X, {\C})]
   \neq 0. $
\end{proposition}

        The following trivial observations were probably  a starting point for Catanese's
         Theorem \ref{GCdF}:

        \begin{lemma}   \label{FibWedge}

(i) For an arbitrary Albanese general compact K\"ahler manifold $ Y $ of
$ \dim _{\C} Y = k, $ the cup product $ \zeta ^{k,0}_{Y} : \wedge ^{k} H^{1,0}(Y)
\rightarrow H^{k}(Y, {\C}) $ is injective and the cup product $ \zeta ^{k+1,0}_{Y} :
\wedge ^{k+1,0} H^{1,0}(Y) \rightarrow H^{k+1}(Y,{\C}) $ is identically zero.

(ii) If $ f_{k} : X \rightarrow Y_{k} $ is an Albanese general $k$-fibration and
$ U_{k} := f^{*}_{k} H^{1,0} (Y_{k}) $ then $ Ker [ \zeta ^{k,0}_{X} : \wedge ^{k}
U_{k} \rightarrow H^{k} (X, {\C})] = 0 $ and $ Im [ \zeta ^{k+1,0}_{X} : \wedge ^{k+1}
U_{k} \rightarrow H^{k+1} (X, {\C}) ] = 0. $ Such subspaces $ U_{k} \subset
 H^{1,0}(X) $ are called strict $k$-wedges.
\end{lemma}

       {\bf Proof: } {\it (i)} According to Proposition \ref{AlbDim},
  there exist $ \omega _{1}, \ldots ,
 \omega _{k} \in H^{1,0}(Y) $ with $ \zeta ^{k,0}_{Y} ( \omega _{1} \wedge  \ldots
 \wedge \omega _{k} ) \neq 0. $ Let $ \{ W ^{(\alpha)} \} _{\alpha \in A} $ be a
 coordinate covering of $ Y. $ For any $ \omega \in H^{1,0}(Y) $ there exist local
 meromorphic functions $ \mu _{i}^{(\alpha)} : W^{(\alpha)} \rightarrow {\C}{\bf
 P}^{1}, $ such that $ \omega \vert _{W^{(\alpha)}} = \sum _{i=1}^{k} \mu
 _{i}^{(\alpha)} \omega _{i} \vert _{W^{(\alpha)}}. $ On the overlaps $ W^{(\alpha)}
 \cap W^{(\beta)} \neq \emptyset, $ one has $ \sum _{i=1}^{k} \left( \mu
 _{i}^{(\alpha)} - \mu _{i}^{(\beta)} \right) \omega _{i} \vert _{W^{(\alpha}) \cap
 W^{(\beta)} } \equiv 0 $ since $ \omega $ and $ \omega _{1}, \ldots , \omega _{k} $
 are globally defined. The linear independence of $ \omega _{1}, \ldots ,
  \omega _{k} $ is inherited by their restrictions on the open subset $ W^{(\alpha)} \cap
  W^{(\beta)} $  of $ Y $ and implies $ \mu _{i}^{(\alpha)} \vert _{ W^{(\alpha)} \cap
   W^{(\beta)}} = \mu _{i}^{(\beta)} \vert _{ W^{(\alpha)} \cap W^{(\beta)} } $ for
   all $ 1 \leq i \leq k. $ In other words, $ \mu _{i} : Y \rightarrow
   {\C}{\bf P}^{1} $ are globally defined and $ \omega \in H^{1,0}(Y) $ can be globally
   represented in the form $ \omega = \sum _{i=1}^{k} \mu _{i} \omega _{i}. $ Consequently,
   $ \wedge ^{k} H^{1,0}(Y) $ consists of $ \mu \omega _{1} \wedge \ldots \wedge
   \omega _{k} $ for some global meromorphic functions $ \mu  : Y \rightarrow
   {\C}{\bf P}^{1}. $ The assumption $ \zeta ^{k,0}_{Y} ( \mu \omega _{1} \wedge
   \ldots \wedge \omega _{k} ) = 0 $ is equivalent to the existence of a
   $(k-1)$-form $ \sigma $ with $ \mu \omega _{1} \wedge \ldots \wedge \omega _{k}
   = d \sigma . $ Then $ 0 = \int _{X} d ( \sigma \wedge d \overline{\sigma}) =
   \int _{X} | \mu |^{2} \omega _{1} \wedge \ldots \wedge \omega _{k} \wedge
   \overline{ \omega _{1}} \wedge \ldots \wedge \overline{\omega _{k}} $ forces the
   vanishing of $ \mu $ almost everywhere on $ Y. $ According to the complex
   analyticity of the zero locus of $ \mu , $ one concludes that $ \mu \equiv 0, $
   i.e., $ \zeta ^{k,0} _{Y} : \wedge ^{k} H^{1,0}(Y) \rightarrow H^{k} (Y,{\C}) $
   is injective. Clearly, $ \zeta ^{k+1,0}_{Y} ( \wedge ^{k+1} H^{1,0} (Y) ) = 0 $ due to
   $ \dim _{\C} Y = k. $

   (ii) The surjective holomorphic map $ f_{k} : X \rightarrow Y_{k} $ induces an
  embedding $ f_{k}^{*}: H^{*}(Y_{k}, {\C}) \rightarrow H^{*}(X, {\C}), $ compatible
  with the cup products. More precisely,
  $$  0 = \zeta ^{l,0}_{X} \left( \sum _{i = (i_{1}, \ldots , i_{l})} f^{*}_{k} ( \omega
  _{i_{1}}) \wedge \ldots \wedge f^{*}_{k}  (\omega _{i_{l}} ) \right) =
  \zeta ^{l,0}_{X} f^{*}_{k} \left( \sum _{i = (i_{1}, \ldots , i_{l})} \omega
  _{i_{1}} \wedge \ldots \wedge \omega _{i_{l}} \right) = $$
$$   f^{*}_{k} \zeta ^{l,0} _{Y_{k}} \left( \sum _{i =( i_{1}, \ldots , i_{l})} \omega
  _{i_{1}} \wedge \ldots \wedge \omega _{i_{l}} \right) $$
  is equivalent to
  $ \zeta ^{l,0}_{Y_{k}} \left( \sum _{i = (i_{1}, \ldots , i_{l})} \omega _{i_{1}}
  \wedge \ldots \wedge \omega _{i_{l}} \right) = 0 $  for any natural number $ l. $
  Putting $ l = k $ or $ k+1 $  and combining with {\it (i)}, one obtains
   {\it  (ii)}, Q.E.D.

      Here is the generalized Castelnuovo-deFranchis Theorem:

   \begin{theorem}    \label{GCdF}
{\rm (Catanese \cite{C})} Let $ X $ be a compact K\"ahler manifold. Then for any
   strict $k$-wedge $ U_{k} \subset H^{1,0}(X) $ there exists an Albanese general
   $k$-fibration $ f_{k} : X \rightarrow Y_{k} $ with $ f^{*}_{k} H^{1,0} (Y_{k}) =
   U_{k}, $ which is unique up to a biholomorphism of $ Y_{k}. $
   \end{theorem}

   \begin{corollary}   \label{AlbGenus}
The $k$-th Albanese genus $ g_{k} $ of a compact K\"ahler manifold $ X $  equals  the maximum
 dimension of a strict $k$-wedge $ U_{k} \subset H^{1,0}(X) $ of $ \dim _{\C} U_{k} > k. $
 \end{corollary}

      In order to formulate one more result of Catanese, used in section 5, let us
     say  that $ V_{k} \subset H^{1,0}(X) $ is a $k$-wedge if
    $ \zeta ^{k,0}_{X} ( \wedge ^{k} V_{k} ) \neq  0 $ and
    $ \zeta ^{k+1,0}_{X} ( \wedge ^{k+1} V_{k}) = 0. $

   \begin{lemma}  \label{Wedge}
{\rm (Catanese \cite{C})} If $ X $ is a compact K\"ahler manifold then any
 $k$-wedge $ V_{k} \subset H^{1,0}(X) $ contains a strict $l$-wedge
 $ U_{l} \subseteq V_{k} $ for some natural number $ l \leq k. $
   \end{lemma}

   \section{Cup products in group and de Rham cohomologies}

   \hspace{.8mm}

        Let us choose a cell decomposition of $ X. $ Then construct an
 Eilenberg-MacLane space $ Y = K( \pi _{1}(X), 1) $ by glueing cells of real
 dimension $ \geq 3 $ to $ X, $ in order to annihilate the higher homotopy groups
 $ \pi_{i}(X), $ $ i \geq 2. $ Put $ c: X \rightarrow Y $ for the resulting
 classifying map and denote by $ S(X)_{\bullet}, $ $ S(Y)_{\bullet} $ the
 corresponding singular chain complexes. The induced chain morphism $ c_{*} :
 S(X)_{\bullet} \rightarrow S(Y)_{\bullet} $ is an isomorphism in degree $ \leq 2 $
 and injective in degree $ i \geq 3. $ If $ \partial _{i}^{X}, $ $ \partial _{i}^{Y} $ are
  the boundary maps $ \partial _{i}^{*} : S(*)_{i} \rightarrow S(*)_{i-1} $
  and $ Z(*)_{i} := \{ \xi \in S(*)_{i} \vert \partial _{i}^{*} ( \xi ) = 0 \} $ are the
  abelian subgroups of the cycles, then $ c_{i} : H_{i}(X,{\Z}) := Z(X)_{i} / \partial _{i+1}^{X}
  S(X)_{i+1} \rightarrow H_{i}(\pi_{1}(X), {\Z}) = H_{i}(Y, {\Z}) := Z(Y)_{i} /
  \partial _{i+1}^{Y} S(Y)_{i+1} $ are isomorphisms for $ i = 0, 1 $ and
   $ c_{2} : H_{2}(X,{\Z}) \rightarrow H_{2}( \pi_{1}(X), {\Z}) $ is surjective. In
   general, the homomorphisms of abelian groups $ c_{i} : H_{i}(X,{\Z}) \rightarrow
   H_{i}( \pi_{1}(X), {\Z}) $ for $ i \geq 3 $ do not obey to any specific
   restrictions.

          For any field $ F $ of $ char F = 0, $ acted trivially by $ \pi_{1}(X), $
the Universal Coefficients Theorems
$$ 0 \rightarrow Ext ^{1}_{\Z} (H_{m-1}(X,{\Z}), F) \rightarrow H^{m}(X,F)
 \rightarrow Hom _{\Z} (H_{m}(X,{\Z}),F) \rightarrow 0, $$
 $$ 0 \rightarrow Ext ^{1}_{\Z} (H_{m-1}( \pi_{1}(X), {\Z}), F) \rightarrow
  H^{m} ( \pi_{1}(X), F) \rightarrow Hom _{\Z} (H_{m}( \pi_{1}(X), {\Z}), F)
  \rightarrow 0, $$
  provide $ H^{m}( \hspace{.5cm}, F) \simeq F^{rk H_{m}( \hspace{.5mm}, {\Z})}, $ due
   to the divisibility of the ${\Z}$-module $ F. $ In particular, the ${\C}$-linear maps
  $ c^{i} : H^{i} ( \pi _{1}(X), {\C}) = H^{i} (Y, {\C}) \rightarrow H^{i}(X,{\C}) $
  are isomorphisms for $ i = 0, 1 $ and injective for $ i = 2. $ Making use of the
  Hodge decomposition $ H^{1} (X, {\C}) = H^{1,0} (X) \oplus H^{0,1} (X) $ on the
  first cohomologies of a compact K\"ahler manifold $ X, $ one introduces
  $ H^{k,l} ( \pi _{1}(X)) := \left( c^{1} \right) ^{-1} H^{k,l}(X) $ for
  $ (k,l) = (1,0) $ or $ (0,1). $ On one hand, there are cup products
  $ \zeta ^{i,j}_{\pi_{1}(X)} : \wedge ^{i} H^{1,0}( \pi_{1}(X)) \otimes _{\C}
  \wedge ^{j} H^{0,1}(\pi _{1}(X)) \rightarrow H^{i+j}(\pi_{1}(X),{\C}) $
  of group cohomologies, defined as a composition of the direct product with the
   dual of a diagonal approximation (cf.\cite{B}).  On the other hand, one has
    cup products $ \zeta ^{i,j}_{X} : \wedge ^{i} H^{1,0} (X) \otimes _{\C}
   \wedge ^{j} H^{0,1} (X) \rightarrow H^{i+j}(X, {\C}) $
 of de Rham cohomologies. Their images are related by the following

 \begin{lemma}    \label{GroupMan}
 Let $ X $ be a compact K\"ahler manifold with fundamental group $ \pi_{1}(X), $
  $  Y = K( \pi_{1}(X), 1) $  be an Eilenberg-MacLane space and
 $ c: X \rightarrow Y $ be a continuous classifying map. Then the cup products
 $$ \zeta ^{i,j}_{\pi_{1}(X)} : \wedge ^{i} H^{1,0}( \pi_{1}(X)) \otimes _{\C}
 \wedge ^{j} H^{0,1}( \pi _{1}(X)) \rightarrow H^{i+j}( \pi_{1}(X), {\C}) \ \ \mbox{  and  } $$
 $$ \zeta ^{i,j}_{X} : \wedge ^{i} H^{1,0}(X) \otimes _{\C}  \wedge ^{j} H^{0,1}(X)
 \rightarrow H^{i+j}(X, {\C}) $$
 have $ c^{i+j}  Im \zeta ^{i,j} _{\pi_{1}(X)} = Im \zeta ^{i,j} _{X}. $ In particular,
   $ rk \zeta ^{i,j} _{\pi_{1}(X)} \geq rk \zeta ^{i,j}_{X} $ for all
  $ i, j \in {\N} \cup \{ 0 \}, $ $ i+j \in {\N}. $
 \end{lemma}

 {\bf Proof:}   As far as  $ Y = K( \pi_{1}(X), 1) $ and $ c : X \rightarrow Y $
 are unique up to homotopy, there is no loss of generality in assuming that $ Y $
 is obtained from $ X $ by glueing cells of real dimension $ \geq 3. $ Then
 $ c: X \rightarrow Y $ and the chain morphism $ c_{*}: S(X)_{\bullet} \rightarrow
 S(Y)_{\bullet} $  are inclusion maps and the dual cochain morphism
  $ c^{*} : S(Y)^{\bullet} :=  Hom _{\Z} ( S(Y)_{\bullet}, {\C}) \rightarrow S(X)^{\bullet} :=
  Hom _{\Z}(S(X)_{\bullet}, {\C}) $ is a surjective restriction from $ Y $ to $ X. $

          By induction on $ i $ one checks that $ c^{i} \zeta ^{i}_{Y} ( \wedge ^{i}
  S(Y)^{1}) = \zeta ^{i}_{X} ( \wedge ^{i} S(X)^{1}) $ for the cup products
  $ \zeta ^{i}_{*} : \wedge ^{i} S(*)^{1} \rightarrow S(*)^{i} $ of cochains. The case
  of $ i = 1 $ is straightforward from the construction of $ Y. $ Any cell
   $ \sigma _{i} \in S(X)_{i} $ is homotopy equivalent to a product of
    segments $ [a_{1}, b_{1}] \times \ldots \times [a_{i}, b_{i}]. $ Let us put
    $ \sigma ' _{i} := [ a_{1}, b_{1}] \times \ldots \times [a_{i-1}, b_{i-1}], $
    $ \sigma ''_{i} := [a_{i}, b_{i}] $ and choose some $ u_{1}, \ldots , u_{i} \in
    S(Y)^{1}. $ Representing the singular cochains on the manifold $ X $ by smooth
    differential forms, one has $ \int _{\sigma ' _{i}} c^{i-1} \zeta ^{i-1}_{\pi_{1}(X)}
     ( u_{1} \wedge \ldots \wedge u_{i-1}) = \int _{ \sigma ' _{i} } \zeta ^{i-1}_{X}
     ( c^{1} (u_{1}) \wedge \ldots \wedge c^{1}(u_{i-1})), $ by the inductional
     hypothesis. Then $ S(X)^{1} = c^{1} S(Y) ^{1} $ and Fubini's Theorem provide
$$ \int _{ \sigma _{i}} \zeta ^{i}_{X} ( c^{1}(u_{1}) \wedge \ldots \wedge
 c^{1}(u_{i}) ) = \int _{\sigma ' _{i} } \zeta ^{i-1}_{X} ( c^{1} (u_{1}) \wedge
 \ldots \wedge c^{1} (u_{i-1})) \int _{\sigma '' _{i}} c^{1}(u_{i}) =  $$
 $$ \int _{\sigma ' _{i}} c^{i-1} \zeta ^{i-1}_{\pi_{1}(X)} ( u_{1} \wedge \ldots \wedge
 u_{i-1}) \int _{\sigma '' _{i}} c^{1} (u_{i}) = \int _{\sigma _{i}} c^{i}
  \zeta ^{i}_{\pi_{1}(X)} ( u_{1} \wedge \ldots \wedge u_{i-1} \wedge u_{i} ) $$
  under the choice of a diagonal approximation $ \Delta : S(Y)_{\bullet} \rightarrow
  S(Y)_{\bullet} \otimes S(Y)_{\bullet}, $ $ \Delta ( [a_{1}, b_{1}] \times \ldots
  \times [a_{m}, b_{m}]) = \sum _{j=0}^{m} ( [a_{1}, b_{1}] \times \ldots \times
  [a_{j}, b_{j}]) \otimes ( [a_{j+1}, b_{j+1}] \times \ldots \times [a_{m}, b_{m}]) $
  on the singular chains of the Eilenberg-MacLane space $ Y = K( \pi_{1}(X), 1). $
  As far as the cup products $ \zeta ^{i}_{X}, $ $ \zeta ^{i}_{Y} $ and the cochain maps
  $ c^{i} : S(Y)^{i} \rightarrow S(X)^{i} $ are ${\C}$-linear, there follows
   $ c^{i} \zeta ^{i}_{Y} ( \wedge ^{i} S(Y)^{1}) = \zeta ^{i}  _{X}
    ( \wedge ^{i} S(X)^{1}) $ for all $ i \in {\N}. $

         The restriction  $ c^{*} : S(Y)^{\bullet} \rightarrow S(X)^{\bullet} $ is a
  morphism of cochain complexes, so that commutes with the coboundary maps
  $ \delta ^{i}_{*} : S(*)^{i} \rightarrow S(*)^{i+1}, $ i.e., $ \delta ^{i}_{X} c^{i}
  = c^{i+1} \delta ^{i}_{Y}. $  In particular, the isomorphisms $ c^{i} : S(Y)^{i}
  \rightarrow S(X)^{i} $ for $ i = 1, 2 $ induce an isomorphism of the 1-cocycles
  $ c^{1} : Z(Y)^{1} \rightarrow Z(X)^{1}, $ where $ Z(*)^{1} := \{ \xi \in S(*)^{1} \ \ \vert
  \ \ \delta ^{1}_{*} ( \xi ) = 0 \}. $ Representing the elements of $ Z(X)^{1} $ by
  closed differential forms and making use of the complex structure $ J $ on $ X, $
  one decomposes $ Z(X)^{1} = Z(X)^{1,0} \oplus Z(X)^{0,1} $ into a direct sum of
  $\pm \sqrt{-1}$-eigenspaces for the action of $ J. $ That allows to introduce
 $ Z(Y)^{k,l} := (c^{1})^{-1} Z(X)^{k,l} $ for $ (k,l) = (1,0) $ or $ (0,1) $ and to
 represent $ Z(Y)^{1} = Z(Y)^{1,0} \oplus Z(Y)^{0,1}. $  Leibnitz' rule for the
 coboundary maps $ \delta ^{i+j}_{X} $ justifies the existence of natural cup products
 $ \zeta ^{i,j}_{X} : \wedge ^{i} Z(X)^{1,0} \otimes _{\C} \wedge ^{j} Z(X)^{0,1}
 \rightarrow Z(X) ^{i+j}. $ As far as $ \wedge ^{i+j} (c^{1})^{-1} : \wedge ^{i} Z(X)
 ^{1,0} \otimes _{\C} \wedge ^{j} Z(X) ^{0,1} \rightarrow \wedge ^{i} Z(Y)^{1,0}
 \otimes _{\C} \wedge ^{j} Z(Y)^{0,1} $ are well defined isomorphisms, compatible with
 $ \delta ^{*}_{X}, $ $ \delta ^{*}_{Y}, $ one can introduce cup products
 $ \zeta ^{i,j}_{Y} : \wedge ^{i} Z(Y) ^{1,0} \otimes _{\C} \wedge ^{j} Z(Y) ^{0,1}
 \rightarrow Z(Y) ^{i+j} $ with   $ c^{i+j} \zeta ^{i,j}_{ Y } ( \wedge ^{i} Z(Y) ^{1,0}
 \otimes _{\C} \wedge ^{j} Z(Y)^{0,1})  =  \zeta ^{i,j} _{X} ( \wedge ^{i} Z(X)^{1,0}
 \otimes _{\C} \wedge ^{j} Z(X) ^{0,1} ). $ On the other hand, the surjectiveness of
 the cochain morphism $ c^{*} $ implies that $ c^{i+j-1} S(Y) ^{i+j-1} = S(X)
 ^{i+j-1}, $ whereas $ c^{i+j} \delta ^{i+j-1} _{Y} S(Y) ^{i+j-1} =
  \delta ^{i+j-1} _{X} c^{i+j-1} S(Y) ^{i+j-1} = \delta ^{i+j-1}_{X}  S(X) ^{i+j-1}. $

   Consequently,
   $$  c^{i+j} Im \zeta ^{i,j}_{ Y } = \frac{c^{i+j} \zeta ^{i,j}_{ Y }
   \left( \wedge ^{i} Z(Y)^{1,0} \otimes _{\C} \wedge ^{j} Z(Y)^{0,1} \right)}
   {c^{i+j} \zeta ^{i,j}_{Y}  \left( \wedge ^{i} Z(Y)^{1,0} \otimes _{\C}
   \wedge ^{j} Z(Y)^{0,1} \right) \cap c^{i+j} \delta ^{i+j} _{Y} S(Y)^{i+j-1} } =  $$
$$ \frac{ \zeta ^{i,j}_{X} \left( \wedge ^{i} Z(X) ^{1,0} \otimes _{\C} \wedge ^{j}
 Z(X)^{0,1} \right) } {  \zeta ^{i,j}_{X}  \left( \wedge ^{i} Z(X) ^{1,0}
  \otimes  _{\C} \wedge ^{j} Z(X) ^{0,1} \right) \cap \delta  ^{i+j-1}_{X} S(X) ^{i+j-1}} =
     Im \zeta ^{i,j} _{X}, \ \ \mbox{  Q.E.D. } $$

 \section{Estimates on cup products}

 \hspace{.8mm}

 The present section provides lower bounds on the rank of cup products of 1-forms on
  $ X.$ Clearly, $ rk \zeta ^{1,0}_{X} = rk \zeta ^{0,1}_{X} = h^{1,0}, $ as far as
$ \zeta ^{i,j}_{X} = Id _{H^{i,j}(X)} $ for $ (i,j) = (1,0) $ or $ (0,1). $

 \begin{lemma}    \label{CupFactorization}
 Let $ X $ be a compact K\"ahler manifold with Albanese dimension $ a. $ Then the
 cup products
 $$ \zeta ^{i,j}_{X} : \wedge ^{i} H^{1,0}(X) \otimes _{\C} \wedge ^{j} H^{0,1}(X)
 \rightarrow H^{i+j} (X,{\C})   $$
 factor through the cup products $ \zeta ^{p,q}_{X} $ for all $ 0 \leq p \leq i , $
 $ 0 \leq q \leq j. $ In particular, $ \zeta ^{i,j}_{X} \equiv 0 $ for $ i > a $ or
 $ j > a $ and the Albanese genera $ g_{k} = 0 $ for all $ k > a. $
 \end{lemma}

 {\bf Proof: } Let us identify the cohomology classes on $ X $ with their de Rham
  representatives. Denote by $ A^{r,s} $ the space of the $C^{\infty}$-forms of type $ (r,s) $
  and put $ Z^{r,s} = \{ \varphi \in A^{r,s} \vert d \varphi = 0 \} $ for the subspace of the
   $d$-closed forms. Then $ \wedge ^{i} H^{1,0}(X) \otimes _{\C} \wedge ^{j} H^{0,1}(X) =
    (\wedge ^{i} Z^{1,0} \otimes  _{\C} \wedge ^{j} Z^{0,1}) / L, $
 $  Im \zeta ^{i,j}_{X} = (\wedge ^{i} Z^{1,0} \otimes _{\C} \wedge ^{j} Z^{0,1}) / M, $
 $ Im  ( \zeta ^{p,q}_{X} \wedge  Id _{\wedge ^{i-p} H^{1,0}(X) \otimes _{\C} \wedge ^{j-q}
 H^{0,1}(X)} =( \wedge ^{i} Z^{1,0} \otimes _{\C} \wedge ^{j} Z^{0,1}) / N $ where
 $$ L := [( d A^{0,0} \cap A^{1,0}) \wedge ( \wedge ^{i-1} Z^{1,0}) ] \otimes  _{\C}
 \wedge ^{j} Z^{0,1}  + \wedge ^{i} Z^{1,0} \otimes _{\C} [ ( d A^{0,0} \cap A^{0,1} ) \wedge
 ( \wedge ^{j-1} Z^{0,1})] ,  $$
 $$ M := ( d A^{i-1,j} + dA^{i,j-1} ) \cap A^{i,j},  \ \
  N := [( d A^{p-1.q} + d A^{p,q-1} ) \cap A^{p,q}] \wedge (  \wedge ^{i-p} Z^{1,0}
 \otimes  _{\C} \wedge ^{j-q} Z^{0,1} ) + $$
  $$ Z^{p,q} \wedge \{ [ ( d A ^{0,0} \cap A^{1,0} ) \wedge ( \wedge ^{i-p-1}Z^{1,0})]
  \otimes   _{\C} \wedge ^{j-q} Z^{0,1} \} + Z^{p,q} \wedge \{  \wedge ^{i-p} Z^{1,0}
   \otimes _{\C} [ ( d A^{0,0} \cap A^{0,1}) \wedge ( \wedge ^{i-q-1} Z^{0,1})] \}. $$
  The existence of correctly defined  ${\C}$-linear maps $ \zeta ^{i,j}_{X}, $
 $ \zeta ^{p,q}_{X} \wedge Id _{ \wedge ^{i-p} H^{1,0}(X) \otimes _{\C} \wedge ^{j-q}
 H^{0,1} (X)} $ for $ 0 \leq p \leq i, $ $ 0 \leq q \leq j $ is due to the inclusions
 $ L  = d ( A^{0,0} \otimes _{\C} \wedge ^{i-1} Z^{1,0} \otimes _{\C} \wedge ^{j} Z^{0,1} )
 \cap A^{i,j} + d ( \wedge ^{i} Z^{1,0} \otimes _{\C} A^{0,0} \otimes _{C} \wedge ^{j-1} Z^{0,1} )
 \cap A^{i,j} \subseteq M, $
  $ L  = ( \wedge ^{p} Z^{1,0} \otimes _{\C} \wedge ^{q} Z^{0,1} ) \wedge \{ [ ( dA^{0,0} \cap
   A^{1,0} ) \wedge ( \wedge ^{i-p-1} Z^{1,0} ) ] \otimes _{\C} \wedge ^{j-q} Z^{0,1} \} +
   ( \wedge ^{p} Z^{1,0} \otimes _{\C} \wedge ^{q} Z^{0,1} ) \wedge \{ \wedge ^{i-p} Z^{1,0}
   \otimes _{\C} [ ( da^{0,0} \cap A^{0,1} ) \wedge ( \wedge ^{j-q-1} Z^{0,1} )] \subseteq N. $
  A necessary and sufficient condition for the factorization of $ \zeta ^{i,j}_{X} $ through
 $ \zeta ^{p,q} _{X} \wedge Id _{\wedge ^{i-p} H^{1,0}(X) \otimes  _{\C}\wedge ^{j-q}
 H^{0,1}(X)} $ is the inclusion
 $ N =  d ( A^{p-1,q} \otimes  _{\C} \wedge ^{i-p} Z^{1,0} \otimes _{\C}  \wedge ^{j-q} Z^{0,1})
 \cap A^{i,j} + d ( A^{p,q-1} \otimes  _{\C} \wedge ^{i-p} Z^{1.0} \otimes _{\C}
  \wedge ^{j-q} Z^{0,1} ) \cap A^{i,j} + d ( Z^{p,q} \otimes  _{\C} A^{0,0} \otimes _{\C}
   \wedge ^{i-p-1} Z^{1,0} \otimes _{\C} \wedge ^{j-q} Z^{0,1} \cap A^{i,j} +
   d ( Z^{p,q} \otimes _{\C} \wedge ^{i-p} Z^{1,0} \otimes _{\C} A^{0,0} \otimes _{\C}
     \wedge ^{j-q-1}Z^{0,1} ) \cap A^{i,j} \subseteq M. $ In particular, $ \zeta ^{i,j}_{X} $
     with $ i > a $ factor through $ \zeta ^{a+1,0}_{X} $ and $ \zeta ^{k,l}_{X} $ with
     $ l > a $ factor through $ \zeta ^{0,a+1}_{X}. $ According to Proposition \ref{AlbDim},
      the cup product $ \zeta ^{a+1,0}_{X} \equiv 0 $ vanishes identically. Hodge duality on the
      compact K\"ahler manifold $ X $ provides $ \zeta ^{0,a+1}_{X} \equiv 0. $ The vanishing of
      $ \zeta ^{k,0}_{X} $ for $ k > a $ implies the nonexistence of strict $k$-wedges
      $ U_{k} \subset H^{1,0}(X) $ (cf. Lemma \ref{FibWedge} {\it (ii)}). Applying Corollary
       \ref{AlbGenus}, one concludes that $ g_{k} = 0 $ for all $ k > a, $  Q.E.D.

 \begin{lemma}     \label{Bounds}
 Let $ X $ be a compact K\"ahler manifold with irregularity $ h^{1,0} > 0, $
  Albanese dimension $ a > 0 $ and  Albanese genera $ g_{k}, $ $ 1 \leq k \leq a. $
   Then the ranks of the cup products
 $$ \zeta ^{i,j} _{X} : \wedge ^{i} H^{1,0} (X) \otimes _{\C} \wedge ^{j} H^{0,1}(X)
 \rightarrow H^{i+j}(X,{\C}) $$
 satisfy the following lower bounds

 (i) $  rk \zeta ^{i,j}_{X} \geq  \left( \begin{array}{c}
                                                   a  \\
                                                   i+j
                                                   \end{array} \right) $
 for  $ i , j \in {\N} \cup \{ 0 \}, $  $  2 \leq i + j \leq a ; $

(ii) $ rk \zeta ^{i,j} _{X} \geq  \left( \begin{array}{c}
                                                   g_{k} - i \\
                                                   j
                                                   \end{array} \right), $
$ rk \zeta ^{j,i}_{X} \geq \left( \begin{array}{c}
                                   g_{k} - i \\
                                   j
                                   \end{array} \right) $
if $ g_{k} > 0 $ for some $ 0 \leq i \leq j, $ $ 2 \leq  i+j \leq k \leq a ; $

(iii) $ rk \zeta ^{i,j}_{X} \geq (i+j) ( h^{1,0} - i - j) + 1 $ if $ g_{k} = 0 $
for all $ 1 \leq k < i + j \leq a; $

(iv) $ rk \zeta ^{1,1}_{X} \geq 2a - 1. $
\end{lemma}

{\bf Proof:} {\it (i)}  By Proposition \ref{AlbDim} one has $ \zeta ^{a,0}_{X} \not \equiv 0. $
 Let $ \omega _{1}, \ldots , \omega _{a} \in  H^{1,0}(X) $ be global holomorphic forms with
  $ \varphi := \zeta ^{a,0}_{X} ( \omega _{1} \wedge \ldots \wedge \omega _{a} ) \neq 0. $
   Consider the subspace
  $$ T^{i,j} := Span _{\C} ( \omega _{t_{1}} \wedge \ldots \wedge \omega _{t_{i}}
  \otimes \overline{ \omega _{s_{1}}} \wedge \ldots \wedge
   \overline{ \omega _{s_{j}}} \ \ \vert \ \  1 \leq t_{1} < \ldots < t_{j} < s_{1}
   < \ldots < s_{j} \leq a ) $$
   of $ \wedge ^{i} H^{1,0}(X) \otimes _{\C} \wedge ^{j} H^{0,1}(X). $ We claim that
 $ T^{i,j} $ is injected by $ \zeta ^{i,j}_{X}, $ so that  $ rk \zeta ^{i,j}_{X} \geq
 \dim _{\C} \zeta ^{i,j}_{X} ( T^{i,j} ) = \dim _{\C} T^{i,j} =
  \left( \begin{array}{c}
   a \\
    i+j
   \end{array}  \right) . $
 Let us observe that  $ \varphi \in H^{a,0}(X) $ is primitive,
 i.e., $ \varphi \wedge \Omega ^{n+1-a} \in H^{n+1,n+1-a}(X) = 0 $ where $ \Omega $
 stands for the K\"ahler form of $ X $ and $ n = \dim _{\C} X. $ By the degeneracy
  of the Hodge Hermitian form on the primitive cohomologies, one has
   $ (-1) ^{\frac{a(a-1)}{2}} \left( \sqrt{-1} \right) ^{a} \int _{X} \varphi \wedge
   \overline{\varphi} \wedge \Omega ^{n-a} > 0 $ and the cohomology class
   $  \zeta ^{a,a}_{X} ( \varphi \wedge \overline{ \varphi })  \neq 0 $ in $ H^{a,a}(X). $
    For any  $ \alpha = \sum _{t,s} c_{t,s} \omega _{t_{1}} \wedge
  \ldots \wedge \omega _{t_{i}} \otimes \overline{ \omega _{s_{1}}} \wedge \ldots
  \wedge \overline{ \omega _{s_{j}}} \in Ker \zeta ^{i,j}_{X} \cap T^{i,j}, $
   $ c_{t,s} \in {\C}, $ let us wedge $ \zeta ^{i,j} _{X} ( \alpha ) =
    d \beta _{1} $ by $ \left( \wedge _{k \in \{ 1, \ldots , a \} \setminus
\{ t_{1}, \ldots , t_{i} \}} \omega _{k} \right) \wedge \left( \wedge _{ l \in
\{ 1, \ldots , a \} \setminus \{ s_{1}, \ldots , s_{j} \}} \overline{ \omega _{l}}
\right), $ to obtain $ \pm c_{t,s} \varphi \wedge \overline{ \varphi} =
 d \beta _{2} $ for appropriate differential forms $ \beta _{1}, \beta _{2}. $ That implies the
 vanishing of all the complex coefficients $ c_{t,s} $ of $ \alpha $ and justifies
 that $ Ker \zeta ^{i,j}_{X} \cap T^{i,j} = 0, $ whereas $ \zeta ^{i,j}_{X} (T^{i,j})
 \simeq T^{i,j}. $

{\it (ii)} If $ g_{k} > 0 $ for some $ i+j \leq k \leq a, $ then according to
Corollary \ref{AlbGenus} there is a strict $k$-wedge $ U_{k} \subset H^{1,0}(X) $
of $ \dim _{\C} U_{k} = g_{k} \geq k+1. $ Let $ u_{1}, \ldots , u_{g_{k}} $ be a
${\C}$-basis of $ U_{k} $ and
 $$ A^{i,j}_{k} := u_{1} \wedge \ldots \wedge u_{i}
\otimes  _{\C}  \wedge ^{j} Span _{\C} \left( \overline{u_{i+1}}, \ldots ,
\overline{ u_{g_{k}}} \right),  \ \  B^{i,j}_{k} := \wedge ^{i}
Span _{\C} \left( u_{j+1}, \ldots , u_{g_{k}} \right) \otimes _{\C} \overline{ u_{1}}
\wedge \ldots \wedge \overline{u_{j}} $$
 be subspaces of $ \wedge ^{i} H^{1,0}(X) \otimes  _{\C} \wedge ^{j} H^{0,1}(X). $
  We claim that $ Ker \zeta ^{i,j}_{X} \cap
A^{i,j}_{k} = 0 $ and $ Ker \zeta ^{i,j}_{X} \cap B^{i,j}_{k} = 0, $ so that
$$ rk  \zeta ^{i,j}_{X} \geq \max \left( \dim _{\C} \zeta ^{i,j}_{X} (A^{i,j}_{k}),
 \dim _{\C} \zeta ^{i,j}_{X} ( B^{i,j}_{k}) \right) =
  \max \left( \dim _{\C} A^{i,j}_{k}, \dim _{\C} B^{i,j}_{k} \right) = $$
$$  \max \left( \left( \begin{array}{c}
        g_{k} - i \\
              j
   \end{array} \right),  \left( \begin{array}{c}
                                    g_{k} - j \\
                                              i
                        \end{array} \right) \right) =   \left( \begin{array}{c}
                                                          g_{k} - \min (i,j)  \\
                                                             \max (i,j)
                                                        \end{array} \right), $$
as far as $ \left( \begin{array}{c}
                     g_{k} - i  \\
                         j
                    \end{array}  \right) : \left( \begin{array}{c}
                                                   g_{k} - j  \\
                                                       i
                                                  \end{array}  \right) =
\prod _{s=1}^{j-i} \left( \frac{g_{k} - j + s  }{i+s} \right) > 1 $ for $ i < j $
and $ g_{k} > i + j. $ Let us assume that $ \psi = u_{1} \wedge \ldots \wedge
u_{i} \otimes _{\C}  \left( \sum _{i+1 \leq s_{1} < \ldots < s_{j} \leq g_{k}} c_{s}
\overline{ u_{s_{1}}} \wedge \ldots \wedge \overline{ u_{s_{j}}} \right) \in
Ker \zeta ^{i,j}_{X} \cap A^{i,j}_{k} $ for some $ c_{s} \in {\C}. $ If $ \Omega $
is the K\"ahler form of $ X $ and $ n = \dim _{\C} X $ then
$ 0 = \int _{X} \psi \wedge \overline{ \psi} \wedge \Omega ^{n-i-j} =
\pm \int _{X} \varphi \wedge \overline{ \varphi} \wedge \Omega ^{n-i-j} $ for
$ \varphi := u_{1} \wedge \ldots \wedge u_{i} \wedge \left( \sum _{i+1 \leq s_{1}
< \ldots < s_{j} \leq g_{k}} c_{s} u_{s_{1}} \wedge \ldots \wedge u_{s_{j}} \right)
\in \wedge ^{i+j} U_{k}. $ As far as $ \zeta ^{n+1,n+1-i-j}_{X} \left( \varphi
\wedge \Omega ^{n+1-i-j} \right) \in H^{n+1,n+1-i-j}(X) = 0, $ the form $ \varphi $
is primitive and $ 0 = \int _{X} \varphi \wedge \overline{\varphi} \wedge
 \Omega ^{n-i-j} = \int _{X} \zeta ^{i+j,0}_X{} ( \varphi ) \wedge \overline{
 \zeta ^{i+j,0}_{X}( \varphi )} \wedge \Omega ^{n-i-j} $ implies that
 $ \zeta ^{i+j,0}_{X} ( \varphi ) = 0, $ according to the nondegeneracy of the Hodge
 Hermitian form on the primitive $ \varphi . $  In other words,
  $ \varphi \in Ker \zeta ^{i+j,0}_{X} \cap \left( \wedge ^{i+j} U_{k} \right). $
   However, $ Ker \zeta ^{k,0}_{X} \cap \left( \wedge ^{k} U_{k} \right) = 0 $
  by the strictness of the $k$-wedge $ U_{k}. $ The factorization of
   $ \zeta ^{k,0}_{X} $ through $ \zeta ^{i+j,0}_{X} $ for $ i + j \leq k $
 implies that $ Ker \zeta ^{i+j,0}_{X} \cap \left( \wedge ^{i+j} U_{k} \right) = 0, $
  whereas $ \varphi = 0. $ Due to the ${\C}$-linear independence
   of $ u_{1} \wedge \ldots \wedge u_{i} \wedge u_{s_{1}} \wedge \ldots \wedge
   u_{s_{j}} $ with $ i+1 \leq s_{1} < \ldots < s_{j} \leq g_{k} $ there follow
   $ c_{s} = 0 $ for all $ s = ( s_{1}, \ldots , s_{j} ). $ Consequently,
   $ Ker \zeta ^{i,j}_{X} \cap A^{i,j}_{k} = 0. $ Similar considerations
  justify $ Ker \zeta ^{i,j}_{X} \cap B^{i,j}_{k} = 0. $

 {\it (iii)} Let us consider the punctured cone
  $$ {\cal C}^{i,j}_{X} := \{ \omega _{1} \wedge \ldots \wedge \omega _{i} \wedge
 \overline{ \omega _{i+1}} \wedge \ldots \wedge \overline{\omega _{i+j}} \ \ \vert
\ \  \omega _{1}, \ldots , \omega _{i+j} \in H^{1,0}(X), \omega _{1} \wedge \ldots
   \wedge \omega _{i} \wedge \omega _{i+1} \wedge \ldots \wedge \omega _{i+j}
    \neq 0 \}. $$
  There is a real diffeomorphism
 $ {\cal C}^{i,j}_{X} \rightarrow {\cal C}^{i+j,0}_{X} := \{ \omega _{1} \wedge
 \ldots \wedge \omega _{i} \wedge \omega _{i+1} \wedge \ldots \wedge \omega _{i+j}
 \neq 0 \vert \omega _{1}, \ldots , \omega _{i+j} \in H^{1,0}(X) \} $
 onto the cone of the decomposable elements of $ \wedge ^{i+j} H^{1,0}(X) $ with
 punctured origin. The projectivization of $ {\cal C}^{i+j,0}_{X} $ is the
 Grassmannian  manifold $ Grass (i+j, H^{1,0}(X)), $ so that $ \dim _{\C} {\cal C}^{i+j,0}_{X} =
 (i+j) (h^{1,0} - i - j) + 1 . $ We claim that $ Ker \zeta ^{i,j}_{X} \cap
 {\cal C}^{i,j}_{X} = \emptyset, $ in order to estimate
 $$ rk \zeta ^{i,j}_{X} = \left( \begin{array}{c}
                                  h^{1,0}   \\
                                  i
                                  \end{array} \right)
\left( \begin{array}{c}
         h^{1,0}  \\
         j
         \end{array}  \right) - \dim _{\C} Ker \zeta ^{i,j}_{X} \geq \dim _{\C}
 {\cal C}^{i,j}_{X} = \dim _{\C} {\cal C}^{i+j,0}_{X}. $$

  Let us assume the opposite,i.e.,   $ \psi = \omega _{1} \wedge \ldots \wedge \omega _{i}
  \wedge \overline{ \omega _{i+1}} \wedge \ldots \wedge \overline{\omega _{i+j}} \in
  Ker \zeta ^{i,j}_{X} \cap {\cal C}^{i,j}_{X}. $ Then
  $ 0 = \int _{X} \psi \wedge \overline{\psi} \wedge \Omega ^{n-i-j} =
  \pm \int _{X} \varphi \wedge \overline{\varphi} \wedge \Omega ^{n-i-j} $
  for $ \varphi = \omega _{1} \wedge \ldots \wedge \omega _{i} \wedge \omega _{i+1}
  \wedge \ldots \wedge \omega _{i+j}. $ On a compact K\"ahler manifold $ X $ of
  $ \dim _{\C} X = n, $ the fact that $ \zeta ^{n+1,n-i-j+1}_{X} ( \varphi \wedge
  \Omega ^{n-i-j+1}) \in  H^{n+1,n-i-j+1}(X) $ reveals the primitiveness of $ \varphi . $ Then the vanishing
  of $ \int _{X} \varphi \wedge \overline{\varphi} \wedge \Omega ^{n-i-j} $ implies  $ \zeta ^{i+j,0}_{X} ( \varphi ) = 0. $ Let $ V_{\varphi } $ be the
   ${\C}$-span of $ \omega _{1}, \ldots , \omega _{i+j}. $ By a decreasing induction
    on   $ 1 \leq k \leq i+j $ will be checked that $ \wedge ^{k} V_{\varphi} \subset
    Ker \zeta^{k,0}_{X}. $ We have already seen that $ {\C} \varphi  = \wedge ^{i+j}
     V_{\varphi} \subset Ker \zeta  ^{i+j,0}_{X}. $ For any natural number
   $ k < i+j $ and a decomposable form $ 0 \neq \psi ' := \omega _{s_{1}} \wedge
   \ldots \wedge \omega _{s_{k}} \in \wedge ^{k} V_{\varphi} , $
   $ \omega _{s_{i}} \in V_{\varphi}, $ there exists
  $ \omega _{s_{k+1}} \in V_{\varphi} $ with $ \psi := \psi ' \wedge
  \omega _{s_{k+1}} \neq 0. $ By the inductional hypothesis $ \psi \in
   \wedge ^{k+1} V_{\varphi} \subset Ker \zeta ^{k+1,0}_{X}. $ If
    $ \zeta ^{k,0}_{X} ( \psi ' ) \neq 0 $ then for $ V_{\psi} := Span _{\C}
    ( \omega _{s_{1}}, \ldots , \omega _{s_{k}}, \omega _{s_{k+1}} ) $ there hold
    $ \zeta ^{k,0}_{X} ( \wedge ^{k} V_{\psi}) \neq 0 $ and $ \zeta ^{k+1,0}_{X}
    ( \wedge ^{k+1} V_{\psi} ) = 0. $ In other words, $ V_{\psi} \subset H^{1,0}(X) $
    appears to be a $k$-wedge and according to Lemma \ref{Wedge}, there is a strict
    $l$-wedge $ U_{l} \subset V_{\psi} $ for some $ 1 \leq l \leq k. $ However,
    $ g_{l} \geq \dim _{\C} U_{l} > 0 $ for $ l < i+j $ contradicts the assumptions
    of {\it (iii)}.  Consequently, $ \zeta ^{k,0}_{X} ( \psi ' ) = 0 $ for all decomposable
     elements of $ \wedge ^{k} V_{\varphi}, $  whereas  $ \wedge ^{k} V_{\varphi} \subset
     Ker \zeta ^{k,0}_{X} $ for all $ 1 \leq k \leq i+j. $ In particular, $ V_{\varphi} \subset Ker \zeta ^{1,0}_{X} = Ker Id
    _{H^{1,0}(X)} $ is an absurd, justifying $ Ker \zeta ^{i,j}_{X} \cap
     {\cal C}^{i,j}_{X} = \emptyset . $ Let us observe that this estimate generalizes
     Amor\'os' Lemma \ref{Amoros1} {\it (i)}.

             {\it (iv)} According to Proposition \ref{AlbDim}, there exist
$ \omega _{1}, \ldots , \omega _{a} \in H^{1,0}(X) $ with $ 0 \neq \omega _{1}
\wedge \ldots \wedge \omega _{a} \wedge \overline{ \omega _{1}} \wedge \ldots
\wedge \overline{\omega _{a}} \in H^{a,a}(X). $   For $ a = 1 $ it is immediate
that $ rk \zeta ^{1,1}_{X} \geq 1. $ For $ a \geq 2 $ we assert that the subspace
 $$ W := Span _{\C}( \omega _{1} \otimes _{\C} \overline{ \omega _{1}}, \ \
  \omega _{1} \otimes _{\C} \overline{ \omega _{i}} ,  \ \ \omega _{i} \otimes _{\C}
   \overline{ \omega _{1}} \ \ \vert \ \ 2 \leq i \leq a) $$
of $ H^{1,0}(X) \otimes _{\C}  H^{0,1}(X)  $ is embedded in $ H^{1,1}(X) $ by the
cup product $ \zeta ^{1,1}_{X}. $ Indeed, if $ \alpha = b_{0} \omega _{1} \wedge
 \overline{\omega _{1}} + \sum _{i=2}^{a} b_{i} \omega _{1} \wedge
  \overline{ \omega _{i}} + \sum _{i=2}^{a} c_{i} \omega _{i} \wedge
\overline{ \omega _{1}} = d \beta $ for some $ b_{0}, b_{i}, c_{i} \in {\C} $
and a 1-form $ \beta, $ then $ \alpha \wedge \overline{\alpha} = d ( \beta \wedge d
\overline{\beta}) = \sum _{i=2}^{a} \sum _{j=2}^{a} ( \overline{b_{i}} b_{j} +
c_{i} \overline{c_{j}} ) \omega _{1} \wedge \overline{\omega _{1}} \wedge
 \omega _{i} \wedge \overline{\omega _{j}}. $ Introducing $ \sigma _{i} :=
 \omega _{2} \wedge \ldots \wedge \omega _{i-1} \wedge \omega _{i+1} \wedge \ldots
 \wedge \omega _{a} $ for $ 2 \leq i \leq a, $ one obtains the vanishing of
 $ \alpha \wedge \overline{\alpha} \wedge \sigma _{i} \wedge \overline{ \sigma _{i}} = d ( \beta
  \wedge d \overline{\beta} \wedge \sigma _{i} \wedge \overline{\sigma _{i}} ) =
 \pm ( |b_{i}|^{2} + |c_{i}|^{2}) \omega _{1} \wedge \ldots \wedge \omega _{a} \wedge
  \overline{ \omega _{1}} \wedge \ldots \wedge \overline{\omega _{a}} \in
 H^{a,a} (X). $ By the choice of $ \omega _{1}, \ldots , \omega _{a}, $ there follow
 $ b_{i} = 0 $ and $ c_{i} = 0 $ for all $ 2 \leq i \leq a, $ whereas $ \alpha =
 b_{0} \omega _{1} \wedge \overline{\omega _{1}} = d \beta . $
  The assumption $ b_{0} \neq 0 $ would imply  $ \omega _{1} \wedge \ldots \omega _{a}
   \wedge \overline{ \omega _{1}} \wedge \ldots \wedge \overline{\omega _{a}}=
   \frac{(-1)^{a-1}}{b_{0}} d \beta \wedge \omega _{2} \wedge \ldots \wedge \omega _{a}
   \wedge \overline{\omega _{2}} \wedge \ldots \wedge \overline{\omega _{a}} =
    d \left( \frac{(-1)^{a-1} \beta}{b_{0}} \wedge \omega _{2} \wedge \ldots \omega _{a}
    \wedge \overline{\omega _{2}} \wedge \ldots \wedge \overline{\omega _{a}} \right)
    = 0 \in H^{a,a}(X),  $ which is an absurd. Therefore $ b_{0} = 0 $ and
     $ Ker \zeta ^{1,1}_{X} \cap W = 0, $ whereas $ rk \zeta ^{1,1}_{X} \geq
      \dim _{\C} \zeta ^{1,1}_{X} (W) = \dim _{\C} W = 2a - 1. $

                As far as    $ \wedge ^{h^{1,0} + 1} H^{1,0}(X) = 0, $ the Albanese dimension
     $ a \leq h^{1,0}. $   Thus, in the case of $ g_{1} = 0 $ Amor\'os' lower bound
       $ rk \zeta ^{1,1}_{X} \geq 2 h^{1,0} - 1 $ from Lemma \ref{Amoros1} {\it (ii)}
      is better than  $ rk   \zeta ^{1,1} _{X} \geq 2a -1, $ Q.E.D.

\section{Proofs of the main results}

\hspace{.8mm}

    The Betti number $ b_{2}(G) := \dim _{\C} H^{2}(G,{\C}) = rk H_{2}(G,{\Z}) $
of an arbitrary finitely presented group $ G $ can be expressed by the means of
 the following

\begin{theorem}     \label{Hopf}
{\rm (Hopf \cite{H}, \cite{B}, \cite{CZ})}  Let $ F =  \langle x_{1} , \ldots , x_{s}
\rangle  $ be a free group, $ R  $  be the normal subgroup of $ F, $ generated by
$ y_{1} , \ldots , y_{r} \in F  $ and $ G = F / R. $  Then there is an exact
 sequence of group homologies
$$ 0 \rightarrow H_{2}(G, {\Z}) \rightarrow H_{1}(R,{\Z})_{G} \rightarrow
 H_{1}(F,{\Z}) \rightarrow H_{1}(G,{\Z}) \rightarrow 0 $$
 where the subscript $ G $ stands for the $G$-coinvariants of the adjoint
 action.
 \end{theorem}

       The ranks of the aforementioned homology groups, will be calculated by the means
  of the following

\begin{lemma}    \label{CoInv}
Let $ F = \langle x_{1}, \ldots , x_{s} \rangle $ be a free group, $ R_{o} $ be the
subgroup of $ F, $ generated by $ y_{1}, \ldots y_{r} \in F, $ $ R $ be the normal
subgroup of $ F, $ generated by $ R_{o} $ and $ G = F / R. $ Then

(i) the coinvariants $ (ab R)_{G} = (ab R)_{F} = R / [F,R] $ are isomorphic to
the image $ R_{o} / ( R_{o} \cap [F,R]) $ of the $F$-coinvariants epimorphism
$ ab R_{o} \rightarrow ab R_{o} / K $ with kernel $ K :=  (R_{o} \cap [F,R]) /
[R_{o},R_{o}]. $ In particular, $ (ab R)_{G} $ is a finitely generated abelian group
of $ rk (ab R)_{G} = r - k $ where $ k := rk K. $

(ii) $ s = rk (ab F) \geq rk (ab (G)) $ with equality exactly for
$ R \subset [F,F]. $
\end{lemma}

{\bf Proof:} Let us recall from \cite{W} the isomorphism $ H_{1}(\Gamma , {\Z})
 \simeq ab \Gamma  := \Gamma / [\Gamma, \Gamma] $ for an arbitrary group $ \Gamma. $
  The adjoint action of $ F $ on its normal subgroup $ R $ descends to an adjoint
  action on $ ab R = R / [R,R], $ as far as $ [R,R] $ is also normal in $ F. $
  Since the adjoint action of $ R $ centralizes $ ab R, $ the $F$-action  on $ ab R $
  coincides with the $G$-action, as well as the corresponding  coinvariants
  $ (ab R)_{F} = (ab R) _{G} . $ The kernel   of the coinvariants epimorphism
  $ ab R \rightarrow (ab R)_{F} $ is generated by  $ f r f^{-1}r^{-1} [R,R] $  for
  $ f \in F, $ $ r \in R. $ Therefore $ (ab R)_{F} = \frac{ab R}{ [F,R] / [R,R] } = R / [F,R]. $

        The normal subgroup $ R $ of $ F $ is generated by $ f^{-1} y_{j} f $ for
$ 1 \leq j \leq r, $ $ f \in F. $  Therefore $ [F,R] $ is generated by
$ (f^{-1}_{1} y_{j} ^{-1} f_{1}) ( f^{-1}_{2} y_{j} f_{2}) $ for $ 1 \leq j \leq r, $
 $ f_{1}, f_{2} \in F. $ In particular, $ y_{j}^{-1} f^{-1} y_{j} f \in [F,R], $ whereas
  $ f^{-1} y_{j} f [F,R] = y_{j}[F,R] $ for all $ 1 \leq j \leq r $ and $ f \in F, $ so that
  any coset  $ r [F,R] \in R/[F,R] $ has a representative $ r_{o} \in  R_{o}, $ $ r_{o} [F,R]
   = r [F,R]. $  That is why, the natural map  $ \psi : R_{o} \rightarrow R / [F,R], $
    $ \psi (r_{o}) =  r_{o} [F,R] $ is an epimorphism with  $ Ker \psi = R_{o} \cap [F,R]. $
     Thus,  $ R_{o} / ( R_{o} \cap [F,R]) $ is isomorphic to $ R / [F,R]. $ Representing
     $ R_{o} / ( R_{o} \cap [F,R]) = (ab R_{o}) / K $ by $ K := (R_{o} \cap [F,R]) /
     [R_{o},R_{o}], $ one concludes that $ R / [F,R] $ is a finitely generated abelian
     group. Clearly, $ rk \left( R / [F,R] \right) =  rk (ab R_{o}) - rk K = r - k. $

  {\it (ii)}    The abelianization is a right exact functor, so that the epimorphism
  $ \alpha : F \rightarrow G $ induces an epimorphism $ \beta : ab F \simeq
  {\Z}^{s} \rightarrow ab G. $ In particular, $ s = rk ( ab F) \geq
  rk ( ab (G)). $ If $ s = rk ( ab (G)) $ then $ \beta $ has to be an
   isomorphism. On one hand, $ F \rightarrow ab F \simeq ab G $ has kernel $ [F,F]. $
   On the other hand, the composed map $ F \rightarrow G  \rightarrow ab G  $
   contains $ R $ in its kernel, so that $ R \subset [F,F]. $ Conversely, if $ R
   \subset [F,F] $ then $ ab G = \frac{G}{[F,F]/R} \simeq ab F, $ Q.E.D.

   {\bf Proof of Proposition \ref{GenRel}: }
           Hopf's Theorem \ref{Hopf} and Lemma \ref{CoInv} imply that
$ b_{2} ( \pi_{1}(X)) = r - k - s + 2 h^{1,0}, $ as far as $ b_{1} ( \pi_{1}(X)) =
b_{1}(X) = 2 h^{1,0}. $  On the other hand, Lemma \ref{GroupMan} provides
$ b_{2}( \pi_{1}(X)) \geq 2 rk \zeta ^{2,0}_{X} + rk \zeta ^{1,1}_{X}. $
The proof is completed by the following immediate consequence of Lemma \ref{Bounds}
and Lemam \ref{Amoros1}{\it (ii)}:

\begin{corollary}      \label{Betti2}
Let $ X $ be a compact K\"ahler manifold  with irregularity $ h^{1,0} > 0,$
Albanese dimension $ a > 0  $ and Albanese genera $ g_{k}, $ $ 1 \leq k \leq a. $
Then the  Betti numbers  $ b_{2} = b_{2} ( \pi_{1}(X)) $ and $ b_{2} = b_{2}(X) =
b_{2n-2}(X) $ are subject to the following lower bounds:

(i) $ b_{2} \geq 1 $ for $ a = 1; $

(ii) $ b_{2} \geq \max \left( a(a-1), g_{k}(g_{k}-1) \ \ \vert \ \ 2 \leq k \leq a
\right) + \max \left( \frac{a(a-1)}{2}, 2a-1, g_{k} - 1 \ \ \vert \ \ 2 \leq k \leq a
\right) $ for $ h^{1,0} \geq g_{1} \geq 2, $ $ a \geq 2; $

(iii) $ b_{2} \geq \max \left( 4 h^{1,0} - 6, a(a-1), g_{k}(g_{k}-1) \ \ \vert \ \
2 \leq k \leq a \right) + $ \\ $ \max \left( 2h^{1,0} - 1, \frac{a(a-1)}{2},
g_{k}-1 \ \ \vert \ \ 2 \leq k \leq a \right) $ for $ h^{1,0} \geq 2, $
$ a \geq 2, $ $ g_{1} = 0. $
\end{corollary}

{\bf Proof of Proposition \ref{BettiNum}: } First of all, let us observe that
$ Im \zeta ^{j,m-j}_{*} \cap \left( \sum _{s \neq j} Im \zeta ^{s,m-s}_{*} \right)
= 0 $ for either of the cup products $ \zeta ^{j,m-j}_{\pi_{1}(X)} $ or
 $ \zeta ^{j,m-j}_{X}. $ Therefore $ b_{m} (*) \geq \sum _{j=0}^{m}
  rk \zeta ^{j,m-j}_{*}. $ Lemma \ref{GroupMan} has established that
  $ rk \zeta ^{j,m-j}_{\pi_{1}(X)} \geq rk \zeta ^{j,m-j}_{X}. $
The lower bounds on $ rk \zeta ^{i,j}_{X} $ from Lemma \ref{Bounds} are invariant
under a permutation of $ i $ with $ j. $ Combining them, one obtains
 $ rk \zeta ^{i,j}_{X} \geq \mu ^{i,j} $ for $ \mu ^{i,j}, $ defined in the statement of
  Proposition \ref{BettiNum},  $ 0 \leq i \leq j, $  $ 3 \leq i+j \leq a. $  Therefore,
  the Betti numbers of the compact K\"ahler manifold $ X $ and its fundamental group
  $ \pi _{1}(X) $ are subject to the inequalities $ b_{2i}(*) \geq 2 \sum _{j=0}^{i-1}
 \mu ^{j,2i-j} + \mu ^{i,i}, $ $ b_{2i+1}(*) \geq 2 \sum _{j=0}^{i} \mu ^{j,2i+1-j} $
 for $ 3 \leq 2i, 2i+1 \leq a. $ In the case of $ 2n - a \leq m \leq 2n - 3, $
 by Serre duality on the cohomologies of $ X $  there hold
  $$ h^{j,m-j}(X) = h^{n-j,n-m+j}(X) \geq  rk \zeta ^{n-j,n-m+j}_{X} \geq
  \mu ^{ n - \max (j, m-j), n - \min (j, n-j)} $$
  as far as $ 3 \leq 2n-m \leq a. $
  Combining with Hodge duality $ h^{n-j, j} (X) = h^{j, n-j}(X) $ for compact
  K\"ahler manifolds $ X, $ one justifies the last two announced inequalities, Q.E.D.

      At first glance, Proposition \ref{GenRel}  can be reformulated entirely in terms
 of the cohomologies of $ \pi _{1}(X). $ However, there are several obstacles for
 doing that. First of all, any isomorphism $ c^{(1)} : H^{1} ( \pi_{1}(X), {\C} )
 \rightarrow H^{1}(X,{\C}) $ allows to introduce $ H^{i,j}( \pi_{1}(X) ) :=
 ( c^{(1)}) ^{-1} H^{i,j}(X) $ for $ (i,j) = (1,0) $ or $ (0,1) $ and to endow
 $ H^{1}( \pi_{1}(X) , {\C}) = H^{1} ( \pi_{1}(X) , {\Z}) \otimes _{\Z} {\C} $ with a
 polarized Hodge structure. For an abstract finitely presented group $ G $ with
 $ H^{1}(G,{\Z}) $ of even rank $ 2q, $ the polarized Hodge structures on
 $ H^{1}(G,{\C}) = H^{1}(G, {\Z}) \otimes _{\Z} {\C} $ are parametrized by the Siegel
 upper half-space $ {\cal S} = Sp(q,{\R}) / U_{q}. $ Unfortunately, the cup products
 $ \zeta ^{i}_{G} : \wedge ^{i} H^{1}(G,{\C}) \rightarrow H^{i} (G,{\C}) $ are not
 invariant under the action of the symplectic group $ Sp(q,{\R}) $ and the
 counterparts of Albanese dimension
 $ a = \max \{ m \in {\N} \cup \{ 0 \} \ \  \vert \ \ \zeta ^{m,0} ( \wedge ^{m}
 H^{1,0}) \neq 0, \zeta ^{m+1,0} ( \wedge ^{m+1} H^{1,0} ) = 0 \} $ and Albanese genera
 $ g_{k} = \max \{ g \in {\N} \cup \{ 0 \} \vert \exists \mbox{  subspace  }
 U \subset H^{1,0}, \dim _{\C} U = g, Ker [ \zeta ^{k} : \wedge ^{k} U \rightarrow
 H^{k} ] = 0, Im [ \zeta ^{k+1} : \wedge ^{k+1} U \rightarrow H^{k+1} ] = 0 \} $
 depend on $ \xi \in {\cal S} $ (cf. Proposition \ref{AlbDim} and Corollary
 \ref{AlbGenus}). The corresponding notions cannot be defined in terms of real cup products.
  Namely, for the real point set $ U^{\R} := Span _{\R} ( u + \overline{u}, \sqrt{-1} u -
   \sqrt{-1} \overline{u} \ \ \vert \ \  u \in U ) $ of $ U \subset H^{1,0}, $ the condition
    $ \zeta ^{2k+1} ( \wedge ^{2k+1} U^{\R} ) = 0 $ is  necessary but not sufficient  for
  $ \zeta ^{k+1} ( \wedge ^{k+1} U ) = 0. $ Finally, the cup product in de Rham
  cohomologies $ H^{*} (X, {\C}) $ of a compact K\"ahler manifold $ X $ with $ \pi_{1}(X) = G $ are
  quotients of the corresponding cup products in $ H^{*}(\pi_{1}(X), {\C}). $ Thus,
  $ \zeta ^{i}_{X} ( \wedge ^{i} H^{1,0} (X)) = 0 $ does not imply
   $ \zeta ^{i}_{\pi_{1}(X)} ( \wedge ^{i} H^{1,0} (\pi_{1}(X))_{s}) = 0, $
    $ \xi \in {\cal S}. $

\vspace{1cm}

\flushleft{
            Section of Algebra \\
            Department of Mathematics and Informatics  \\
            Kliment Ohridski University of Sofia  \\
            5 James Bouchier Blvd., Sofia 1126 \\
            e-mail:kasparia@fmi.uni-sofia.bg(Kasparian)}


\newpage
\begin{thebibliography}{99}
\bibitem{A}      Amor\'os J., On the Malcev completion of K\"ahler groups,
                 Comment. Math. Helv. {\bf 71} (1996), 192-212.
\bibitem{BPV}    Barth W.,C. Peters, A. Van de Ven , {\it Compact Complex
                 Surfaces}, Springer Verlag, New York - Heidelberg - Berlin,
                 1984.
\bibitem{ABCKT}  Amor\'os J., M. Burger, K. Corlette, D. Kotschick, D. Toledo,
                 {\it Fundamental Groups of Compact K\"ahler Manifolds,} Math.
                 Surv. Mon. {\bf 44,} Providence, 1996.
\bibitem{B}      Brown K.S., {\it Cohomology of Groups}, Springer Verlag,
                 New York - Heidelberg - Berlin, 1982.
\bibitem{C}      Catanese F., Moduli and classification of irregular K\"ahler
                 manifolds (and algebraic varieties) with Albanese general type
                 fibrations, Inv. Math. {\bf 104} (1991), 263-289.
\bibitem{CZ}     Collins D., H. Zieschang, Combinatorial group theory and
                 fundamental groups, {\it Encyclopedia of Mathematics,} vol.
                 {\bf 58} (1990), 5-100 (in Russian).
\bibitem{GL}     Green M., R. Lazarsfeld, Higher obstructions to deforming
                 cohomology groups of line bundles, Jour. AMS {\bf 4} (1991)
                 87-103.
\bibitem{H}      Hopf H., Fundamentalgruppe und Zweite Bettische Gruppe, Comment.
                 Math. Helv. {\bf 14} (1942) 257-309.

\bibitem{RV}     Remmert R., A. Van de Ven, Zur Funktionentheorie homogener
                  komplexer Mannigfaltigkeiten, Topology {\bf 2} (1963),
                  137-157.
\bibitem{W}      Weibel C., {\it An Introduction to Homological Algebra, }
                 Cambridge University Press, Cambridge, 1995.
\end{thebibliography}
\end{document}